\documentclass [11pt]{article}
\usepackage[english]{babel}
\usepackage[latin1]{inputenc} %escribe con acentos y eñes
\usepackage{amssymb,amsmath}

\newtheorem{theorem}{Theorem}[section]
\newtheorem{corollary}[theorem]{Corollary}
\newtheorem{lemma}[theorem]{Lemma}
\newtheorem{proposition}[theorem]{Proposition}
\newtheorem{remark}{Remark}
\newtheorem{definition}[theorem]{Definition}

\usepackage{hyperref}

\title{The canonical involution in the space of connections of a $(J^{2}=\pm 1)$-metric manifold}
\author{Fernando Etayo\footnote{Dept. Mathematics, Statistics and Computation.  University of Cantabria. Avda. de los Castros, s/n, 39071 Santander, SPAIN. e-mail: etayof@unican.es} \, and Rafael Santamar\'{\i}a\footnote{Departamento de Matem\'{a}ticas. Escuela de Ingenier\'{\i}as Industrial e Inform\'{a}tica. Universidad de Le\'{o}n. Campus de Vegazana, 24071 Le\'{o}n, SPAIN. e-mail: rsans@unileon.es}}
\date{}
\begin{document}
\maketitle

\begin{abstract}
A $(J^{2}=\pm 1)$-metric manifold has an almost complex or almost product structure $J$ and a compatible metric $g$. 
We show that there exists  a canonical involution in the set of connections on such a manifold, which allows to define a projection over the set of connections adapted to $J$. This projection sends the Levi Civita connection onto the first canonical connection. In the almost Hermitian case, it also sends the $\nabla^{-}$ connection onto the Chern connection, thus applying the line of metric connections defined by $\nabla ^{-}$ and the Levi Civita connections onto the line of canonical connections. Besides, it moves metric connections onto metric connections.
\end{abstract}



\section{Introduction}

In the celebrated paper \cite{gauduchon} of Gauduchon, connections on almost Hermitian manifolds are studied, focusing on the $1$-parameter family of canonical connections, which is defined as 
\begin{equation}
\nabla^{\mathrm{t}}= (1-t) \nabla^{0} + t \nabla^{\mathrm{c}},   \quad \forall t \in \mathbb{R},
\label{eq:conexionescanonicas}
\end{equation}
where $\nabla^{0}$ denotes the first canonical connection and $\nabla^{\mathrm{c}}$ the Chern connection. It is also said that $\nabla^{0}$ is the orthogonal projection of the Levi Civita connection $\nabla^{\mathrm{g}}$ onto the affine space of Hermitian connections. As is well known,  both connections coincide in the K\"{a}hler case.

The purpose of this note is threefold: (1) we want to extend the above result to all the $(J^{2}=\pm 1)$-metric manifolds, i.e., manifolds endowed with an almost complex or almost product structure $J$ and a compatible metric $g$,  (2) we will show, in the almost Hermitian context, that the connection $\nabla ^{-}$ projects onto the Chern connection, and (3) we will prove that the projection moves metric connections onto metric connections. 
\bigskip

Connections with totally skew-symmetric torsion are very useful in Physics (see, e.g., \cite{agricola, friedrich} and the references therein). In particular, connections $\nabla^{+}$ and $\nabla ^{-}$ appear in heterotic string theory (see, e.g \cite{heterotic, nuclear} and the references therein).  Having a non-K\"{a}hler manifold of type $\mathcal G_1$ in the classification of almost Hermitian manifolds of Gray and Hervella \cite{gray-hervella}, there exists a unique Hermitian connection $\nabla^{\mathrm{sk}}$ with totally skew-symmetric torsion (see \cite[Theor. 5.21]{brno}). Then one can define the connections
\begin{align*}
\nabla^{+}_X Y& = \nabla^{\mathrm{g}}_X Y + \frac{1}{2} \mathrm{T}^{\mathrm{sk}}(X,Y), \quad \forall X, Y \in \mathfrak{X} (M),\\
\nabla^{-}_X Y &= \nabla^{\mathrm{g}}_X Y - \frac{1}{2} \mathrm{T}^{\mathrm{sk}}(X,Y), \quad \forall X, Y \in \mathfrak{X} (M),
\end{align*}
where $\mathrm{T}^{\mathrm{sk}}$ denotes the torsion tensor of $\nabla^{\mathrm{sk}}$. We will show that $\nabla ^{+}$ is invariant under the projection and that $\nabla ^{-}$ projects onto the Chern connection.
\bigskip

We will consider smooth manifolds and operators being of class $C^{\infty}$. As in this introduction, $\mathfrak{X}(M)$ denotes the module of vector fields of a manifold $M$.

\section{Canonical connections}

We are dealing with all the four geometries: almost Hermitian, almost Norden, almost product Riemannian and almost para-Hermitian, which correspond to the cases 
\[
(\alpha ,\varepsilon )\in\{(-1,1), (-1,-1),\ (1,1), (1,-1)\}
\]
in the following

\begin{definition}[{\cite[Defin. 3.1]{racsam}}]
Let $M$ be a  manifold, $g$ a semi-Riemannian metric on $M$, $J$ a  tensor field of type (1,1) and $\alpha ,\varepsilon \in \{-1,1\}$. Then $(J,g)$ is called an $(\alpha ,\varepsilon )$-structure on $M$ if
\[
J^{2} = \alpha Id, \quad \mathrm{trace} \, J=0, \quad g(JX,JY)= \varepsilon g(X,Y), \quad \forall X, Y \in \mathfrak{X}(M), 
\]
$g$ being a Riemannianan metric if $\varepsilon =1$. Then  $(M,J,g)$ is called a $(J^{2}=\pm 1)$-metric manifold.
\end{definition}

Condition $\mathrm{trace} \, J=0$ is a consequence of the other conditions in all the cases unless the $(1,1)$. We impose it in this case looking for a common treatment of all the four geometric structures. See  \cite{racsam} for a more complete description. 
\bigskip

A linear connection $\nabla$ is said to be adapted to the metric $g$ (resp. to $J$) if $\nabla g=0$ (resp. $\nabla J=0$). We use the following notation: $\mathcal{C}(M)$ (resp. $\mathcal{C}(M,J)$, $\mathcal{C}(M,g)$, $\mathcal{C}(M,J,g)=\mathcal{C}(M,J)\cap \mathcal{C}(M,g)$) denotes the affine space of linear connections on $M$ (resp. adapted to $J$, to $g$ and to both $J$ and $g$). In \cite{brno} we have studied a lot of distinguished connections defined in such a manifold. Then $\nabla ^{\mathrm{g}}\in \mathcal{C}(M,J,g)$ if and only if $(M,J,g)$ is of K\"{a}hler type. We are interested in the non K\"{a}hler type case. Then two adapted connections will be essential in our study: the first canonical connection $\nabla ^{0}$ and the Chern connection $\nabla^{\mathrm{c}}$. 

Connection $\nabla ^{0}$ has been introduced in \cite{brno} as
\begin{definition}
\label{teor:first-connection}
Let $(M,J,g)$ be a $(J^2=\pm 1)$-metric manifold. The first canonical connection of $(M,J,g)$ is the linear connection having the covariant derivative $\nabla ^{0}$ given by
\[
\nabla^{0}_X Y= \nabla^{\mathrm{g}}_X Y +\frac{(-\alpha)}{2} (\nabla^{\mathrm{g}}_X J)JY,   \quad \forall X,Y \in \mathfrak{X} (M).
\]
\end{definition}
The previous one generalizes the classical definition given in the context of almost Hermitian manifolds (see, e.g., \cite{gauduchon}). 

The Chern connection was firstly introduced in the case of Hermitian manifolds \cite{chern}. In  \cite{racsam} we have extended the connection to the almost para-Hermitian case, recovering the connection defined by Cruceanu and one of us in \cite{etayo}. The following results establish the existence and uniqueness of the Chern connection on a $(J^2=\pm1)$-metric manifold with $\alpha\varepsilon=-1$.

\begin{theorem} [{\cite[Theor. 6.3]{racsam}}]
\label{teor:chern-connection}
Let  $(M,J,g)$ be a $(J^2=\pm1)$-metric manifold with  $\alpha\varepsilon =-1$. Then there exists a unique linear connection on $M$ adapted to $(J,g)$  defined by $(J,g)$  whose torsion tensor  $\mathrm{T}^{\mathrm{c}}$ satisfies the following condition
\[
\mathrm{T}^{\mathrm{c}}(JX,JY)= \alpha \mathrm{T}^{\mathrm{c}}(X,Y),   \quad \forall X, Y \in \mathfrak{X} (M).
\]
This connection is called the Chern connection of $(M,J,g)$.
\end{theorem}

\begin{remark}

The Chern connection can not be defined in the $\alpha \varepsilon =1$ context, as we have proved in \cite[Remark 6.4]{racsam}. In \cite{brno} we have proved that  in the case $\alpha \varepsilon =-1$, the so-called well adapted connection $\nabla^{\mathrm{w}}$ is also a canonical connection, i.e., is a connection in the line defined in (\ref{eq:conexionescanonicas}). Then this line can be parametrized as 
$\nabla^{\mathrm{s}}= (1-s) \nabla^{0} + s \nabla^{\mathrm{w}},  \forall s \in \mathbb{R}$. As the first canonical connection and the well adapted connection can be also defined in the case $\alpha \varepsilon =1$, we have been able to define canonical connections on any $(J^{2}=\pm1)$-metric manifold $(M,J,g)$. 

\end{remark}

\begin{remark} \label{coincide} (1) We have seen in \cite[Remark 6.2]{brno}, assuming  $\alpha\varepsilon =-1$, that the Chern connection corresponds to the case $s=3$ and in \cite[Example 6.3]{brno} that the Bismut connection $\nabla ^{\mathrm{b}}$ (see \cite{bismut}) to $s=-3$. This connection coincides, if there exists, with the unique adapted connection with totally skew-symmetric torsion. In the almost Hermitian case the well adapted connection ($s=1$) coincides with the connection of minimal torsion defined by Gauduchon  in \cite{gauduchon}.

(2) In the case $\alpha\varepsilon=1$, if $J$ is a non-integrable $\alpha$-structure (which is equivalent  to $\nabla^0\neq\nabla^{\mathrm{w}}$ according to \cite[Theor. 5.6]{brno}), and $(M, J, g)$ is a quasi-K\"{a}hler type manifold, then there exists a unique canonical connection with totally skew-symmetric torsion,  which is that given by $s = -1$.
\end{remark}

\section{Canonical involution and projection of connections}

First of all, we will need the following well known results of Affine Geometry.

\begin{itemize}
\item A subset of an affine space is an affine subspace if and only if the line joining any pair of points of the subset is contained in the subset.

\item A map between affine spaces is an affine map if and only if it preserves barycentric combinations.

\item An involutive affine map $\sigma$ in an affine space defines a projection $\pi =\frac{1}{2}Id +\frac{1}{2}\sigma$ onto the subspace of fixed points of $\sigma$. The map $\pi$ is also an affine map.
\end{itemize}

Taking into account the above properties one easily checks that $\mathcal{C}(M)$, $\mathcal{C}(M,J)$, $\mathcal{C}(M,g)$, $\mathcal{C}(M,J,g)$ and the $1$-parameter family of canonical connections are affine spaces. We introduce the following

\begin{definition} Let $M$ be a manifold endowed with a tensor field $J$ of type $(1,1)$ such that $J^{2}=\alpha Id$, where $\alpha \in \{ -1,1\}$. The  map $J_{*}: \mathcal{C} (M)\to \mathcal{C}(M)$ defined as 
\[
(J_{*}(\nabla ))_{X}(Y)= \alpha J(\nabla _{X}(JY)), \quad \forall X, Y \in \mathfrak{X} (M),
\]
is called the canonical involution induced by $J$ in the affine space of connections  $\mathcal{C} (M)$.
\end{definition}

Then we have:

\begin{proposition} Let $M$ be a manifold endowed with a tensor field $J$ of type $(1,1)$ such that $J^{2}=\alpha Id$, where $\alpha \in \{ -1,1\}$.

\begin{enumerate}
\item The map $J_{*}: \mathcal{C} (M)\to \mathcal{C} (M)$ is an involutive affine isomorphism.

\item ${\mathcal C} (M,J)=\{ \nabla \in \mathcal{C} (M): J_{*}(\nabla )=\nabla \}$.

\item For each $\nabla \in \mathcal{C} (M)$, the connection
\[
\pi (\nabla )= \frac{1}{2} \nabla + \frac{1}{2} J_{*}(\nabla )
\]
is $J$-invariant, so defining a projection $\pi :\mathcal{C} (M)\to \mathcal{C} (M,J) $.

\item In addition, if $(M,J,g)$ is a $(J^{2}=\pm 1)$-metric manifold, then $\pi (\nabla ^{\mathrm{g}})=\nabla ^{0}$.

\end{enumerate}

\end{proposition}

\textit{Proof.}

\begin{enumerate}

\item Direct calculations show that $J_{*}(\nabla)$ is a covariant derivative and
\[
J_{*}(\lambda _{1}\nabla ^{(1)}+\ldots +\lambda _{k}\nabla ^{(k)}) =\lambda _{1}J_{*}(\nabla ^{(1)})+\ldots +\lambda _{k}J_{*}(\nabla ^{(k)})
\]
when $\lambda _{1}+\ldots +\lambda _{k}=1$, thus proving $J_{*}$ is an affine map.  Besides, given $X, Y$ vector fields on $M$ one has:
\[(J_{*}(J_{*}(\nabla )))_{X}(Y)=\alpha J((J_{*}(\nabla ))_{X}(JY))= \alpha J(\alpha J(\nabla _{X} (JJY)))= \alpha ^{4}\; \nabla _{X}Y=\nabla _{X}Y,
\]
thus proving $J_{*}$ is an involutive affine isomorphism.

\item Let $\nabla \in \mathcal{C} (M,J)$. Then
\[
(J_{*}(\nabla ))_{X}(Y)= \alpha J(\nabla _{X}(JY)) =
\alpha J^{2}(\nabla _{X}Y)= \alpha ^{2}\; \nabla _{X}Y = \nabla _{X}Y, \quad \forall X, Y \in \mathfrak{X}(M).
\]
The reverse: Suppose that  $J_{*}(\nabla )=\nabla$ then 
\[
\nabla _{X}Y=(J_{*}(\nabla ))_{X}(Y)= \alpha J(\nabla _{X}(JY)), \quad \forall X, Y \in \mathfrak{X}(M), 
\]
therefore 
\[
J(\nabla _{X}Y) = J(\alpha J(\nabla _{X}(JY)))= \alpha J^{2}(\nabla _{X}(JY))= \nabla _{X}(JY), \quad \forall X, Y \in \mathfrak{X}(M).
\]

\item Taking into account the above items one has:
\[
J_{*}(\pi (\nabla ))= J_{*}\left(\frac{1}{2} \nabla + \frac{1}{2} J_{*}(\nabla )\right)= \frac{1}{2}  J_{*}(\nabla )+ \frac{1}{2}  J_{*}(J_{*}(\nabla )) = \pi (\nabla ).
\]

\item Let $X, Y$ be vector fields on $M$, a straightforward calculation is enough:
\begin{align*}
 \nabla ^{0}_{X}Y &=\nabla ^{\mathrm{g}} _{X}Y +\frac{(-\alpha )}{2}(\nabla ^{\mathrm{g}}_{X}J)JY = \nabla ^{\mathrm{g}}_{X}Y +\frac{(-\alpha )}{2}(\nabla ^{\mathrm{g}}_{X}(J^{2}Y) - J(\nabla ^{\mathrm{g}}_{X}(JY)))\\
 &= \nabla ^{\mathrm{g}}_{X}Y-\frac{1}{2}\nabla ^{\mathrm{g}}_{X}Y + \frac{1}{2}\alpha J(\nabla ^{\mathrm{g}}_{X}(JY))=  \frac{1}{2} \nabla ^{\mathrm{g}} _{X}Y+ \frac{1}{2} (J_{*}(\nabla ^{\mathrm{g}}))_{X}Y =(\pi (\nabla ^{\mathrm{g}}))_{X}Y.  \quad \Box
\end{align*}
\end{enumerate}

\begin{remark} In the case of having a K\"{a}hler type manifold $(M,J,g)$ then  $\nabla ^{0}=\nabla ^{\mathrm{g}}$. In the non-K\"{a}hler type case, $\nabla ^{0}$ is an adapted connection to $(J,g)$ which is obtained as  the projection of $\nabla ^{\mathrm{g}}$ to the set $\mathcal{C}(M,J)$. In fact, $\nabla ^{0}\in \mathcal{C}(M,J,g)$.

\end{remark}

\begin{corollary} Let $(M,J,g)$ be an almost Hermitian non K\"{a}hler manifold of type $\mathcal G_1$. Then 
\[
\pi (\nabla ^{+})=\nabla ^{+}, \quad \pi (\nabla ^{-})=\nabla ^{\mathrm{c}}.
\]
\end{corollary}

\textit{Proof.} Taking into account Remark \ref{coincide},  $\nabla^{+}=\nabla^{\mathrm{sk}}=\nabla^{\mathrm{b}}$ is the Bismut connection, which belongs to the line of canonical connections. As this line is contained in ${\cal C}(M,J)$, one obtains $\pi (\nabla ^{+})=\nabla ^{+}$.

Observe that $\nabla^{\mathrm{g}}=\frac{1}{2}\nabla ^{+} + \frac{1}{2}\nabla ^{-}$. Then 
\[
\nabla ^{0} = \pi (\nabla^{\mathrm{g}})=\frac{1}{2}\pi (\nabla ^{+}) + \frac{1}{2}\pi (\nabla ^{-})= \frac{1}{2}\nabla ^{+} + \frac{1}{2}\pi (\nabla ^{-})= \frac{1}{2}\nabla ^{\mathrm{b}} + \frac{1}{2}\pi (\nabla ^{-}).
\]As $\nabla ^{0}$ is the midpoint between $\nabla^{\mathrm{b}}$ and 
$\nabla^{\mathrm{c}}$, one obtains $\pi (\nabla ^{-})=\nabla ^{\mathrm{c}}$. $\Box$
\bigskip

In \cite{nuclear} the authors have studied the plane of connections defined by the line of canonical connections on an almost Hermitian non K\"{a}hler type manifold and the Levi Civita connection. This plane has another significant line defined by $\nabla^{\mathrm{g}}, \nabla ^{+} , \nabla ^{-}$, when there exists a connection with totally skew-symmetric torsion. We have seen that the projection $\pi$ applies this line of connections onto the line of canonical connections.

\section{Metric connections}

Let $(M,J,g)$ be a $(J^{2}=\pm 1)$-metric manifold and let $\nabla$ be a metric connection. We are going to prove that $\pi (\nabla )$ is also a metric connection. First of all, we obtain the following new expression of the projection $\pi$.

\begin{lemma} Let $(M,J,g)$ be a $(J^{2}=\pm 1)$-metric manifold. The projection $\pi : \mathcal C (M) \rightarrow  \mathcal C (M,J) $ is given by $\pi (\nabla) = \nabla + S_{\nabla}$, where 
\[
S_{\nabla} (X, Y)= \frac{(-\alpha)}{2} (\nabla_X J)JY, \quad \forall X, Y \in \mathfrak{X}(M).
\]
\end{lemma}

\textit{Proof.}  A direct calculus shows that:
\begin{align*}
\nabla_X Y + S_{\nabla}(X, Y)&= \nabla_X Y + \frac{(-\alpha)}{2} (\nabla_X J)JY = \nabla_X Y -\frac{1}{2} \nabla_X Y + \frac{\alpha}{2} (J \nabla_X (JY)) \\
&= \frac{1}{2} \nabla_X Y + \frac{\alpha}{2} (J \nabla_X (JY)) =
 (\pi(\nabla))_X Y, \quad \forall X, Y \in \mathfrak{X}(M). \quad \Box
 \end{align*}
 \smallskip
 
 We need another lemma. 
 \begin{lemma} Let $(M,J,g)$ be a $(J^{2}=\pm 1)$-metric manifold, $\nabla$  a metric connection and let $S$ be  a $(1,2)$ tensor field. Then the connection $\nabla +S$ is  metric  if and only if 
\[
g(S(X,Y),Z)+g(S(X,Z),Y)=0, \quad \forall X, Y, Z \in \mathfrak{X}(M).
\]
\end{lemma}

\textit{Proof.} Given $X, Y, Z$ vector fields on $M$, as $\nabla g=0$, one has: 
\[
X(g(Y,Z))=g(\nabla_X Y, Z)+ g(\nabla_X Z, Y).
\] 
Then $((\nabla + S)_X g)(Y,Z)=0$ if and only if
\[
X(g(Y,Z))=g(\nabla_X Y, Z) + g(S(X,Y),Z)+ g(\nabla_X Z, Y)+g(S(X,Z),Y),
\]
and thus one can conclude the result. $\Box$
 
\begin{proposition} 
Let $(M,J,g)$ be a $(J^{2}=\pm 1)$-metric manifold and $\nabla$  a metric connection. Then $\pi (\nabla )$ is also a metric connection.
 \end{proposition}
 
 \textit{Proof.} In order to apply the above lemma, given $X, Y, Z$ vector fields on $M$ we obtain:
 \begin{align*}
g(S_{\nabla}(X,Y),Z)&= \frac{(-\alpha)}{2} g((\nabla_X J)JY,Z) \\
&= \frac{(-\alpha)}{2} (\alpha g(\nabla_X Y, Z)-g(J(\nabla_X(JY)),Z)) \\
&= \frac{(-\alpha)}{2} (\alpha g(\nabla_X Y, Z)-\alpha \varepsilon g(\nabla_X(JY),JZ))\\
g(S_{\nabla}(X,Z),Y)&= \frac{(-\alpha)}{2} (\alpha g(\nabla_X Z, Y)-\alpha \varepsilon g(\nabla_X(JZ),JY)),
\end{align*}
and then
\begin{align*}
g(S_{\nabla}(X,Y),Z)+g(S_{\nabla}(X,Z),Y)&=- \frac{1}{2}(g(\nabla_X Y, Z)+g(\nabla_X Z, X))\\
&+ \frac{\varepsilon}{2} (g(\nabla_X(JY),JZ) + g(\nabla_X(JZ),JY))\\
&= \frac{1}{2}X(g(Y,Z))+\frac{\varepsilon}{2} X(g(JY,JZ)) \\
&= -\frac{1}{2}X(g(Y,Z))+\frac{1}{2}X(g(Y,Z)) =0. \quad \Box
\end{align*}

\medskip

Thus,  $\pi (\mathcal C(M,g)) \subset \mathcal C(M,g)\cap \mathcal C(M,J) = \mathcal C(M, J,g)$, moving metric connections onto  connections adapted to $(J,g)$. In the case of the plane of connections considered in \cite{nuclear} in the almost Hermitian context, the plane remains globally invariant under the projection $\pi$, moving all the points to the line of canonical connections. More in general, $\pi$ moves metric connections onto Hermitian connections.

\bigskip

\textbf{Acknowledgments.} The authors are grateful to their colleagues L. Ugarte and R. Villacampa for their useful comments.


\begin{thebibliography}{99}

\bibitem{agricola} I.\ Agricola.\ The Srn\'{\i} lectures on non-integrable geometries with torsion.\  {\em Archivum Math.\ Brno.\ }  {\bf 42} (2006), 5--84.


\bibitem{bismut}  J.\ M.\ Bismut.\  A local index theorem for non K\"{a}hler manifolds.\ {\em Math.\ Ann.\ }  {\bf 284} (1989), 681--699.


\bibitem{chern} S.\ S.\ Chern.\ Characteristic classes of Hermitian manifolds.\ {\em Ann.\ of Math.\ } {\bf 47} (1946), 85--121.

\bibitem{etayo} V.\  Cruceanu and  F.\  Etayo.\  On almost para-Hermitian manifolds.\    {\em Algebras Groups Geom.\ }  {\bf 16} no.\ 1 (1999) 47--61.

\bibitem{brno} F. Etayo and R. Santamar\'{\i}a.  Distinguished connections on $(J^{2}=\pm 1)$-metric manifolds.\ {\em Archivum Math.\ Brno.\ } {\bf 52} no.\ 3 (2016), 159--203. 

\bibitem{racsam} F. Etayo and R. Santamar\'{\i}a.  The well adapted connection of a $(J^2=\pm1)$-metric manifold.\ {\em RACSAM.\ } {\bf 111} no.\ 2 (2017), 355--375.

\bibitem{heterotic} M. Fern\'{a}ndez, S. Ivanov, L. Ugarte and D. Vassilev. Non-Kaehler heterotic string solutions with non-zero fluxes and non-constant dilaton.  {\em J. High Energy Phys. }  (2014), 2014: 73.

\bibitem{friedrich} T.\ Friedrich and S.\ Ivanov.\ Parallel spinors and connections with skew-symmetric torsion in string theory.\ {\em Asian J.\ Math.\ }  {\bf 6} (2) (2002),  303--335.

\bibitem{gauduchon} P.\ Gauduchon.\ Hermitian connections and Dirac operators.\  {\em Boll.\ Un.\ Mat.\ Ital.\ B } {\bf 7} suppl. 11 (2) (1997), 257--288.

\bibitem{gray-hervella} A.\ Gray and L.\ M.\ Hervella.\  The sixteen classes of almost Hermitian manifolds and their linear invariants.\ \begin{em} Ann.\ Mat.\ Pura  Appl.\  \end{em}{\bf 123}  no.\ 1 (1980), 35--58.


\bibitem{nuclear} A. Otal, L. Ugarte and R. Villacampa.\ Invariant solutions to the Strominger system and the heterotic equations of motion.  \textit{Nuclear Phys. B} \textbf{920} (2017), 442--474.

\end{thebibliography}
\end{document}